\documentclass[11pt]{article}
\usepackage[usenames]{color}
\usepackage{amssymb,psfig,epsfig,latexsym,graphicx}
\usepackage{amsmath}

\definecolor{webgreen}{rgb}{0,.5,0}
\definecolor{webbrown}{rgb}{.6,0,0}

\newtheorem{theorem}{Theorem}

\newcommand{\eqn}[1]{(\ref{#1})}
\newcommand{\bsq}{{\vrule height .9ex width .8ex depth -.1ex }}

\newcommand{\GG}{{\mathnormal G}}
\newcommand{\HH}{{\mathnormal H}}

\newcommand{\ga}{\gamma}

\newcommand{\et}{\eta}
\newcommand{\si}{\sigma}
\newcommand{\sG}{{\mathcal G}}
\newcommand{\eeq}{\end{equation}}

\newcommand{\beql}[1]{\begin{equation}\label{#1}}

\makeatletter
\def\@sect#1#2#3#4#5#6[#7]#8{\ifnum #2>\c@secnumdepth
     \def\@svsec{}\else
     \refstepcounter{#1}\edef\@svsec{\csname the#1\endcsname.\hskip .75em }\fi
     \@tempskipa #5\relax
      \ifdim \@tempskipa>\z@
        \begingroup #6\relax
          \@hangfrom{\hskip #3\relax\@svsec}{\interlinepenalty \@M #8\par}%
        \endgroup
       \csname #1mark\endcsname{#7}\addcontentsline
         {toc}{#1}{\ifnum #2>\c@secnumdepth \else
                      \protect\numberline{\csname the#1\endcsname}\fi
                    #7}\else
        \def\@svsechd{#6\hskip #3\@svsec #8\csname #1mark\endcsname
                      {#7}\addcontentsline
                           {toc}{#1}{\ifnum #2>\c@secnumdepth \else
                             \protect\numberline{\csname the#1\endcsname}\fi
                       #7}}\fi
     \@xsect{#5}}
\def\@begintheorem#1#2{\it \trivlist \item[\hskip \labelsep{\bf #1\ #2.}]}

\def\section{\@startsection {section}{1}{\z@}{-3.5ex plus -1ex minus
 -.2ex}{2.3ex plus .2ex}{\normalsize\bf}}
\makeatother
\makeatletter
\def\subsection{\@startsection {subsection}{1}{\z@}{-3.5ex plus -1ex minus
 -.2ex}{2.3ex plus .2ex}{\normalsize\bf}}

\makeatother

\begin{document}


\begin{center}
{\large\bf The Gift Exchange Problem } \\
\vspace*{+.2in}

David Applegate and N. J. A. Sloane${}^{(a)}$, \\
AT\&T Shannon Labs, \\
180 Park Ave., Florham Park, \\
NJ 07932-0971, USA.

\vspace*{+.2in}
${}^{(a)}$ Corresponding author.
\vspace*{+.2in}

Email: david@research.att.com, njas@research.att.com. \\
\vspace*{+.2in}
July 1, 2009
\vspace*{+.2in}


{\bf Abstract}
\end{center}

The aim of this paper is to solve the ``gift exchange'' problem:
you are one of $n$ players, and there are 
$n$ wrapped gifts on display; when your turn comes,
you can either choose any of the remaining wrapped gifts,
or you can ``steal'' a gift from someone who has already unwrapped it, subject
to the restriction that no gift can be stolen more than a total
of $\si$ times.  The problem is to determine the number of ways
that the game can be played out, for given values of $\si$ and $n$.
Several recurrences and explicit formulas are given for these
numbers, although some open questions remain.

\vspace{0.8\baselineskip}

\vspace{0.8\baselineskip}


\section{The problem}\label{Sec1}

The following game is sometimes played at parties.
A number $\si$ (typically $1$ or $2$) is fixed in advance.
Each of the $n$ guests brings a wrapped gift, the gifts are
placed on a table (this is the ``pool'' of gifts),
and slips of paper containing the numbers $1$ to $n$
are distributed randomly among the guests.
The host calls out the numbers $1$ through $n$ in order.

When the number you have been given is called, you can either choose one 
of the wrapped (and so unknown) gifts remaining
in the pool, or you can take (or ``steal'') a gift that some earlier
person has unwrapped, subject to the restriction that no gift can be
``stolen'' more than a total of $\si$ times.

If you choose a gift from the pool, you unwrap it and show it
to everyone.
If a person's gift is stolen from them, they immediately get
another turn, and can either take a gift from the pool, or can
steal someone else's gift, subject always to the limit
of $\si$ thefts per gift.
The game ends when someone takes the last ($n$th) gift.

The problem is to determine the number of possible ways
the game can be played out,
for given values of $\si$ and $n$.

For example, if $\si=1$ and $n=3$, with guests $A, B, C$ and 
gifts numbered $1$, $2$, $3$, there are 42 different scenarios,
as follows.  We write $XN$ to indicate that 
guest $X$ took gift $N$ -- it is always clear from the context
whether the gift was stolen or taken from the pool. 
Also, provided we multiply the final answer by 6, we
can assume that the gifts are taken from the pool in 
the order $1, 2, 3$.
There are then seven possibilities:
\begin{align}\label{Eq1}
 {} & A1, B2, C3 \nonumber \\
 {} & A1, B2, C1, A3 \nonumber \\
 {} & A1, B2, C1, A2, B3 \nonumber \\
 {} & A1, B2, C2, B3 \nonumber \\
 {} & A1, B2, C2, B1, A3 \nonumber \\
 {} & A1, B1, A2, C3 \nonumber \\
 {} & A1, B1, A2, C2, A3 \nonumber \\
\end{align}
and so the final answer is $6 \cdot 7 = 42$.

If we continue to ignore the factor of $n!$ due to
the order in which the gifts are selected from the pool,
the number of scenarios for the case $\si=1$ and $n=1,2,3,4,5$ are
$1,2,7,37,266$, respectively.

We noticed that these five terms matched the beginning of entry
A001515 in \cite{OEIS}, although indexed differently.
The $n$th term of A001515 is defined as $y_n(1)$,
where $y_n(x)$ is a Bessel polynomial
(\cite{Gross78}, \cite{KF49}, \cite{RCI}),
and for $n=0,1,2,3,4$ the values are $1,2,7,37,266$, respectively.
Although there was no mention of gift-swapping in that entry,
one of the comments there provided
enough of a hint to lead us to a complete solution
of the general problem.

\subsection*{Comments on the rules}

(i) If $\si = 1$ then once a gift has been stolen it
can never be stolen again.

\noindent
(ii) If $\si = 2$, and someone steals your gift, 
then if you wish you may immediately steal it back (provided
you got it honestly!), and then
it cannot be stolen again. Retrieving a gift in this way, 
although permitted by a strict interpretation of the rules,
may be prohibited at real parties.

\noindent
(iii) A variation of the game allows the last player
to take {\em any} gift that has been unwrapped,
regardless of how many times it has already been stolen,
as an alternative to taking the last gift from the pool.
This case only requires minor modifications of
the analysis, and we will not consider it here.

\noindent
(iv) We also ignore the complications caused by the fact that 
you brought (and wrapped) one of the gifts yourself, and so are
presumably unlikely to choose it when your number is called.

\section{Connection with partitions of labeled sets}\label{Sec2}

Let $\HH_{\si}(n)$ be the number of scenarios with $n$ gifts and a limit
of $\si$ steals, for $\si \ge 0, n \ge 1$.
Then $\HH_{\si}(n)$ is a multiple of $n!$, and we write
$\HH_{\si}(n) = n! \GG_{\si}(n-1)$, where in $\GG_{\si}(n-1)$ we
assume that the gifts are taken from the pool in the order $1,2,\ldots,n$.
We write $n-1$ rather than $n$ as the argument of $\GG_{\si}$ because the
$n$th gift plays a special (and less important) role. This also
simplifies the statement of Theorem \ref{th1}.

In other words, $\GG_{\si}(n)$ is the number of scenarios when there
are $n+1$ gifts, with a limit of $\si$ steals per gift, and 
the gifts are taken from the pool in the order $1,2,\ldots,n+1$. 

As mentioned above, the sequence of values of $\GG_1(n)$ appeared to
coincide with entry A001515 in \cite{OEIS}.
One of the interpretations of that sequence (contributed by
Robert A. Proctor on April 18, 2005) involved partitions of a 
labeled set into blocks, and this was enough of a hint to lead us to
our first theorem.

We recall that the Stirling number of the second kind,
$S_2(i,j)$, is the number of partitions of the labeled set $\{1,\ldots,i\}$ into $j$ blocks
(\cite{Comtet}, \cite{GKP}),
while for $h \ge 1$ the $h$-restricted Stirling number of the second kind,
$S_2^{(h)}(i,j)$, is the number of partitions of $\{1,\ldots,i\}$ into $j$ blocks
of size at most $h$ (\cite{ChSm1}-\cite{ChSm3}).

\begin{theorem}\label{th1}
For $\si \ge 0$ and $n \ge 0$,
\beql{Eq2}
\GG_{\si}(n) = \sum_{k=n}^{(\si+1)n} S_2^{(\si+1)}(k,n)\,.
\eeq
\end{theorem}

\noindent{\bf Proof.}
Equation \eqn{Eq2} is an assertion about $\GG_{\si}(n)$, so we are
now discussing scenarios where there are $n+1$ gifts.
For $\si = 0$, $\HH_0(n+1) = (n+1)!$, so $\GG_0(n) = 1$,
in agreement with $S_2^{(1)}(n,n) = 1$.

We may assume therefore that $\si \ge 1$.
Let an ``action'' refer to a player choosing a gift $\ga$, either by
taking it from the pool or by stealing it from another player.  Since
we are now assuming that the gifts are taken from the pool in order,
$\ga$ determines both the player and whether the action was to take
a gift from the pool or to steal it from another player.  So the
scenario is fully specified
simply by the sequence of $\ga$ values,
recording which gift is chosen at each action.
For example, the scenarios in \eqn{Eq1} are represented
by the sequences
$123$, $1213$, $12123$, $1223$, $12213$, $1123$, $11223$.
Since the game ends as soon as the $(n+1)$st gift is selected, the number
of actions is at least $n+1$ and at most $(\si+1)n+1$.

The sequence of $\ga$ values is therefore a sequence of integers
from $\{1,\ldots,n+1\}$ which begins with $1$, ends with $n+1$,
where each number $i \in \{1,\ldots,n\}$ appears 
at least once and at most $\si+1$ times
and $n+1$ appears just once,
and in which the first $i$ can appear only after $i-1$ has appeared.
Conversely, any sequence with these
properties determines a unique scenario.

Let $k$ denote the length of the sequence with the last entry
(the unique $n+1$) deleted.
We map this shortened sequence to a partition of $[1,\ldots,k]$
into $n$ blocks: the first block records the positions of the $1$'s,
the second block records the positions of the $2$'s, $\ldots$,
and the $n$th block records the positions of the $n$'s.
Continuing the example,
for the seven sequences above,
the values of $k$ and the corresponding partitions are as
shown in Table 1.

\begin{table}[htb]
\label{Tab1}
\caption{Values of $k$ and partitions corresponding to
the scenarios in \eqn{Eq1}. }
$$
\begin{array}{cc}
k & \mbox{partition} \\
2 & 1, 2 \\
3 & 13, 2 \\
4 & 13, 24 \\
3 & 1, 23 \\
4 & 14, 23 \\
3 & 12, 3 \\
4 & 12, 34 \\
\end{array}
$$
\end{table}

The number of such partitions is precisely $S_2^{(\si+1)}(k,n)$.
Since the mapping from sequences to partitions is completely reversible,
the desired result follows.~~~$\bsq$

\vspace*{+.2in}
\noindent{\bf Remark.} The sums $B(i) := \sum_{j}^{} S_2(i,j)$ 
are the classical Bell numbers. The sums $\sum_{j}^{} S_2^{(h)}(i,j)$
also have a long history \cite{MMW}, \cite{MW55}.
However, the sums
$\sum_{i}^{} S_2^{(h)}(i,j)$
mentioned in \eqn{Eq2} do not seem to
have studied before.
Note that the limits in \eqn{Eq2}
are the natural limits on the summand $k$, and
could be omitted.

To simplify the notation, and to put the most important variable first, let
\beql{EqE}
E_{\si}(n,k) :=  S_2^{(\si+1)}(k,n)\,,
\eeq
for $\si \ge 0$, $n \ge 0$, $k \ge 0$.
In words, $E_{\si}(n,k)$ is the number of partitions of $\{1, \ldots, k\}$ into exactly
$n$ blocks of sizes in the range $[1, \ldots, \si+1]$.

For $n \ge 0$, $E_{\si}(n,k)$ is nonzero only
for $n \le k \le (\si+1)n$.
To avoid having to worry about negative arguments,
we define $E_{\si}(n,k)$ to be zero if either $n$ or $k$ is negative.
Then
\beql{Eq4}
\GG_{\si}(n) = \sum_{k=n}^{(\si+1)n} E_{\si}(n,k)\,.
\eeq

Stirling numbers of the second kind satisfy many different recurrences
and generating functions (\cite[Chap.~V]{Comtet}),
and to a lesser extent this is also true for $E_{\si}(n,k)$.
We begin with three general properties.

\begin{theorem}\label{th2}

(i) Suppose $\si \ge 1$. Then $E_{\si}(n,k) = 0$ for $k<n$ or $k>(\si+1)n$,
and otherwise, for $n \le k \le (\si + 1)n$,
\beql{EqAAB}
E_{\si}(n,k)
= \sum_{i=0}^{\si} \binom{k-1}{i} E_{\si}(n-1,k-1-i) \,.
\eeq

(ii) For $\si \ge 0$, $n \ge 0$, $k \ge 0$,
\beql{EqAAA}
E_{\si}(n,k)
= \sum_{(a_1,\ldots,a_{\si+1})}
\frac{k!}{a_1! a_2! \ldots a_{\si+1}! \, 1!^{a_1} 2!^{a_2} \cdots (\si+1)!^{a_{\si+1}} } \,,
\eeq
where the sum is over all $(\si+1)$-tuples of 
nonnegative integers $(a_1,\ldots,a_{\si+1})$ satisfying
\begin{eqnarray}\label{EqAAD}
a_1 + a_2 + a_3 \cdots + a_{\si+1} & = & n \,, \nonumber \\
a_1 + 2a_2 + 3a_3 \cdots + (\si+1)a_{\si+1} & = & k \,. 
\end{eqnarray}

(iii)
The numbers $E_{\si}(n,k)$ have the exponential generating function
\beql{EqAAC}
\sum_{n=0}^{\infty} \sum_{k=n}^{(\si+1)n} E_{\si}(n,k)x^n \frac{y^k}{k!}
=
\exp \left[ x\left(y+\frac{y^2}{2!}+\cdots + \frac{y^{\si+1}}{(\si+1)!}\right) \right] \,.
\eeq
\end{theorem}

\noindent{\bf Proof.}
(i) This is an analog of the ``vertical'' recurrence
for the Stirling numbers
(\cite[Eq.~$\lbrack$3c$\rbrack$,~p.~209]{Comtet}).
The idea of the proof is to take a partition of $[1,\ldots,k]$,
remove the block containing $k$, and renumber the
remaining parts.
(ii) Here $a_i$ is the number of blocks of size $i$ in the partition.
This follows by standard counting 
arguments (cf. \cite[Th.~B,~p.~205]{Comtet}).
(iii) This is an analog of the ``vertical'' generating function
for the Stirling numbers (\cite[Eq.~$\lbrack$2b$\rbrack$,~p.~206]{Comtet}),
and follows directly from (i).~~~$\bsq$

\vspace*{+.2in}
The recurrence in Theorem \ref{th2}(i) makes it easy to  
compute as many values of $E_{\si}(n,k)$ as one wishes.
Tables 3 through 7 show the initial values
of $E_{1}(n,k)$ through $E_{5}(n,k)$, and
Table 8 gives the
initial values of $\GG_{\si}(n)$ for $\si =0$ through $8$.

\section{The case $\si = 1$}\label{Sec3}

In the case when a gift can be stolen at most once, from Theorem \ref{th2}
we have the recurrence
\beql{EqE1a}
E_1(n,k) = E_1(n-1,k-1) + (k-1)E_1(n-1,k-2) \,,
\eeq
for $n \le k \le 2n$, with $E_1(n,k)=0$ for $k<n$ and $k>2n$;
the explicit formula
\beql{EqE1b}
E_1(n,k) = \frac{k!}{(2n-k)!~(k-n)!~2^{k-n}} \,,
\eeq
for $n \le k \le 2n$; and the generating function
\beql{EqE1c}
\sum_{n=0}^{\infty} \sum_{k=n}^{2n}~E_{1}(n,k)~x^n \frac{y^k}{k!}
=
e^{x(y+y^2/2)} \,.
\eeq
It follows from \eqn{Eq4} that
\begin{eqnarray}\label{EqG1a}
\GG_1(n) & = & \sum_{k=n}^{2n} \frac{k!}{(2n-k)!~(k-n)!~2^{k-n}} \nonumber \\
       & = & \sum_{i=0}^{n} \frac{(n+i)!}{(n-i)!~i!~2^i}  \,.
\end{eqnarray}

Equation \eqn{EqG1a} shows that the sequence $\GG_1(n)$ is indeed given by 
entry A001515 in \cite{OEIS}.
That entry gives (mostly without proof) several other properties of these numbers,
taken from various sources, notably Grosswald \cite{Gross78}.
We collect some of these properties in the next theorem.
Property (iii) is especially interesting, since the following sections
will be concerned with attempts to generalize it to larger values of $\si$.
We recall from \cite{Gross78} that the Bessel polynomial $y_n(z)$ is given
by
\beql{EqBess1}
y_n(z) := \sum_{i=0}^{n} \frac{(n+i)!z^i}{(n-i)!~i!~2^i}  \,.
\eeq
Also ${}_2F_{0}$ and (later) ${}_2F_{1}$ denote hypergeometric functions.
\begin{theorem}\label{th3}
(i) 
\beql{EqG1b}
\GG_1(n) = y_n(1)\,.
\eeq
(ii) 
\beql{EqG1c}
\GG_1(n) = {}_2F_{0}\left[ \begin{array}{c}
                            n+1,-n \\
                             -
                            \end{array}
                             ;
                            \begin{array}{c}
                            -\frac{1}{2}
                            \end{array}
                            \right] \,.
\eeq
(iii)
\beql{EqG1d}
\GG_1(n) = (2n-1)\GG_1(n-1) + \GG_1(n-2) \,.
\eeq
for $n \ge 2$, with $\GG_1(0)=1, \GG_1(1)=2$.

\noindent
(iv)
\beql{EqG1e}
\sum_{n=0}^{\infty} \GG_1(n)\frac{x^n}{n!} ~=~ \frac{ e^{1-\sqrt{1-2x}}}{\sqrt{1-2x}} \,.
\eeq
(v)
\beql{EqG1f}
\GG_1(n) ~ \sim ~ \frac{e(2n)!}{n! 2^n} \mbox{~as~} n \rightarrow \infty \,. 
\eeq
\end{theorem}

\noindent{\bf Proof.}
(i) and (ii) are immediate consequences of \eqn{EqG1a}. 

\noindent
(iii) We give three proofs of \eqn{EqG1d}.
(First proof.) Equation \eqn{EqG1d} follows from one of the recurrences for Bessel
polynomials (\cite[Eq.~(7),~p.~18]{Gross78}, \cite{KF49}).
(Second proof.) Alternatively, it is easy to verify from \eqn{EqE1b} that  
\beql{EqE1d}
E_1(n,k) = (2n-1)E_1(n-1,k-2) + E_1(n-2,k-2)\,.
\eeq
Our conventions about negative arguments make it 
unnecessary to put any restrictions on the range over which \eqn{EqE1d}
holds. By summing \eqn{EqE1d} on $k$ we obtain \eqn{EqG1d}.
(Third proof.) The third proof is combinatorial. We will show the equivalent
statement that for $n \ge 3$,
\beql{EqG1g}
\GG_1(n) = \GG_1(n-2) + \GG_1(n-1) + 2(n-1)\GG_1(n-1)\,.
\eeq
We can build a partition counted in $\GG_1(n)$ in three ways.
(A) Take a partition $P$ into $n-2$ parts
and adjoin two parts of size $1$, $\{x\}$ and $\{y\}$, say, where
$x$, $y$ are elements not in $P$.
This gives $\GG_1(n-2)$ partitions.
(B) Take a partition $P$ into $n-1$ parts
and adjoin a part $\{x,y\}$ of size $2$.
This gives $\GG_1(n-1)$ partitions.
(C) Let $P$ be a partition into $n-1$ parts
and let $S$ be one of the parts.
If $S = \{u\}$ is a singleton, then
$$
P \setminus S \cup \{u,x\} \cup \{y\} \mbox{~and~}
P \setminus S \cup \{u,y\} \cup \{x\}
$$
are two partitions into $n$ parts.
If $S = \{u,v\}$ is a pair, then
$$
P \setminus S \cup \{u,x\} \cup \{v,y\} \mbox{~and~}
P \setminus S \cup \{u,y\} \cup \{v,x\}
$$
are two partitions into $n$ parts.
So in either case the pair $P$, $S$ gives rise to two
partitions into $n$ parts.
There are $n-1$ choices for $S$, so in all we obtain $2(n-1)\GG_1(n-1)$
partitions. 
The argument is clearly reversible, and so \eqn{EqG1g} and hence \eqn{EqG1d}
follow.

\noindent
(iv) Let 
\begin{eqnarray}\label{EqG1h}
\sG _1(x) & := & \sum_{n=0}^{\infty} \GG_1(n) \frac{x^n}{n!} \nonumber \\
& = & 1 +2x + 7\frac{x^2}{2!} + 37\frac{x^3}{3!} + 266 \frac{x^4}{4!} + \cdots \,.  \nonumber 
\end{eqnarray}
By multiplying \eqn{EqG1d} by $x^n/n!$ and summing on $n$ from $2$
to $\infty$ we obtain the differential equation
\beql{EqG1i}
\sG_1''(x) = 3 \sG_1'(x) + 2x\sG_1''(x) + \sG_1(x)\,.
\eeq
Then the right-hand side of \eqn{EqG1e} is the unique solution of \eqn{EqG1i}
which satisfies $\sG_1(0) = 1$, $\sG_1'(0) = 2$.

\noindent
(v) This follows from \eqn{EqG1a}, since the terms $i=n-1$
and $i=n$ dominate the sum (see also \cite[Eq.~(1),~p.~124]{Gross78}).~~~$\bsq$

\section{The case $\si = 2$}\label{Sec4}

In the case when a gift can be stolen at most once, the problem, as we saw in the
previous section, turned out to be related to the values of Bessel polynomials,
and the principal sequence, $\GG_1(n)$, had been studied before.
For $\si \ge 2$, we appear to be in new territory---for one thing,
the sequences $\GG_2(n), \GG_3(n), \ldots$ were not among the 140,000 existing
sequences in \cite{OEIS}.

These sequences can be computed using Theorem \ref{th2}.
From \eqn{Eq4}, \eqn{EqAAA} we have:
\beql{EqGsa}
\GG_{\si}(n) ~=~ \sum_{k=n}^{(\si+1)n} \sum_{(a_1,\ldots,a_{\si+1})}
\frac{k!}{a_1! a_2! \ldots a_{\si+1}! \, 1!^{a_1} 2!^{a_2} \cdots (\si+1)!^{a_{\si+1}} } \,,
\eeq
where the inner sum is over all $(\si+1)$-tuples of 
nonnegative integers $(a_1,\ldots,a_{\si+1})$ satisfying \eqn{EqAAD}.
This may be rewritten as a sum of multinomial coefficients:
\beql{EqGsb}
\GG_{\si}(n) ~=~
\frac{1}{n!}~
\sum_{i_1=1}^{\si+1}
\sum_{i_2=1}^{\si+1}
\cdots
\sum_{i_n=1}^{\si+1}
\genfrac{(}{)}{0pt}{0}{i_1+i_2+\cdots+i_{n}}{i_1,~i_2,~\cdots,~i_{n}} \,,
\eeq
where $i_r$ is the size of the $r$th part.

We naturally tried to find analogs of the various parts of Theorem \ref{th3}
that would hold for $\si \ge 2$.
Let us begin with the simplest result, the asymptotic behavior.
This is directly analogous to Theorem \ref{th3}(v).

\begin{theorem}\label{th4}
For fixed $\si \ge 1$,
\beql{EqGsf}
\GG_{\si}(n) ~ \sim ~  \frac{e((\si+1)n)!}{n! {(\si+1)!}^n} \mbox{~as~} n \rightarrow \infty \,.
\eeq
\end{theorem}

\noindent{\bf Sketch of proof.}
The two terms corresponding to
$ \{ k=(\si+1)n, a_{\si+1}=n$, other $a_i=0 \} $ and
$ \{ k=(\si+1)n-1, a_{\si+1}=n-1, a_{\si}=1$, other $a_i=0 \} $
dominate the right-hand side of \eqn{EqGsa},
and are both equal to 
$((\si+1)n)!/(n! {(\si+1)!}^n)$.
Dividing the sum by this quantity gives a converging sum,
in which a subset of terms approach $1+1+1/2!+1/3!+...$,
while the others vanish as $n \rightarrow \infty$.~~~$\bsq$

Concerning Theorem \ref{th3}(i), we do not know if there is a generalization 
of Bessel polynomials whose value gives \eqn{EqGsa} for $\si \ge 2$.

As for Theorem \ref{th3}(ii), there is a relationship with hypergeometric
functions in the case $\si = 2$.
From \eqn{EqAAA} we have
\begin{eqnarray}\label{EqE2a}
E_2(n,k) & = & 
\sum_{c=\max\{0,k-2n\}}^{\lfloor(k-n)/2\rfloor}
\frac{k!}
{(2n-k+c)! (k-n-2c)! c! \, 2^{k-n-c} 3^c} \nonumber \\
& = &
\sum_{c=\max\{0,\et-n\}}^{\lfloor \et/2\rfloor}
\frac{k!}
{(n-\et+c)! (\et-2c)! c! \, 2^{\et-c} 3^c} \,,
\end{eqnarray}
where $\et=k-n$ (this is the ``excess'' of $k$ over $n$).

\begin{theorem}\label{th5}
(i) Let $\et=k-n$.

\noindent
If $\et \le n$ then
\beql{EqE2b}
E_2(n,k) =
\frac{(n+\et)!}{\et! (n-\et)! 2^{\et} } ~ 
   {}_2F_{1}\left[\begin{array}{c}
                  -\et/2,-\et/2+1/2 \\
                       n-\et+1
                  \end{array}
                   ;
                  \begin{array}{c}
                  \frac{8}{3}
                  \end{array}
                  \right] \,.
\eeq

\noindent
If $\et \ge n$ then
\beql{EqE2c}
E_2(n,k) =
\frac{(\et+n)!}{(2n-\et)! (\et-n)! 2^{n} 3^{\et-n} } ~ 
  {}_2F_{1}\left[ \begin{array}{c}
                  -n+\et/2,-n+\et/2+1/2 \\
                       \et-n+1
                  \end{array}
                   ;
                  \begin{array}{c}
                  \frac{8}{3}
                  \end{array}
                  \right] \,.
\eeq

\noindent
(ii)
\begin{eqnarray}\label{EqG2c}
\GG_2(n) & = &
\sum_{ \et = 0 }^{n-1} ~
\frac{(n+\et)!}{\et! (n-\et)! 2^{\et} } ~ 
  {}_2F_{1}\left[ \begin{array}{c}
                  -\et/2,-\et/2+1/2 \\
                       n-\et+1
                  \end{array}
                   ;
                  \begin{array}{c}
                  \frac{8}{3}
                  \end{array}
                  \right]  \nonumber \\
& + & \sum_{ \et = n }^{2n} ~
\frac{(n+\et)!}{(2n-\et)! (\et-n)! 2^{n} 3^{\et-n} } ~ 
  {}_2F_{1}\left[ \begin{array}{c}
                  -n+\et/2,-n+\et/2+1/2 \\
                       \et-n+1
                  \end{array}
                   ;
                  \begin{array}{c}
                  \frac{8}{3}
                  \end{array}
                  \right] \,.
\end{eqnarray}
\end{theorem}

\noindent{\bf Proof.}
(i) follows from \eqn{EqE2a} using the standard rules for converting sums of
products of factorials to hypergeometric functions (cf. \cite{And74}),
and (ii) follows from \eqn{Eq4}.~~~$\bsq$

We can now state the main theorem of this section, which gives
analogs of \eqn{EqE1d} and \eqn{EqG1d}.

\begin{theorem}\label{th6}
(i)
\begin{align}\label{EqE2d}
 E_2(n,k) &  = (9 n^2 - 9 n + 2) E_2(n-1,k-3)/2
       - 5 E_2(n-1,k-1)/2 \nonumber \\
& +\, (9 n^2 - 36 n + 35) E_2(n-2,k-4)/2
       + 6(n-1) E_2(n-2,k-3)
       - 3 E_2(n-2,k-2)/2 \nonumber \\
& +\, 3(2 n-5) E_2(n-3,k-4)
       + 5 E_2(n-3,k-3)/2
       + 5 E_2(n-4,k-4)/2 \, .
\end{align}

\noindent
(ii) 
\begin{align}\label{EqG2d}
\GG_2(n) & = (9 n^2 - 9 n - 3) \GG_2(n-1)/2 \nonumber \\
& +\, (9 n^2 - 24 n + 20) \GG_2(n-2)/2 \nonumber \\
& +\, (6 n - 25/2) \GG_2(n - 3) + 5 \GG_2(n - 4)/2 \, ,
\end{align}
for $n \ge 4$, with $\GG_2(0)=1$, $\GG_2(1) = 3$,
$\GG_2(2) = 31$, $\GG_2(3) = 18252$.
\end{theorem}

\noindent{\bf Proof.}
(ii) Eq. \eqn{EqG2d} follows by summing \eqn{EqE2d}
on $k$, just as \eqn{EqG1d} followed from \eqn{EqE1d}.

(i) We give two proofs of \eqn{EqE2d}.
The first proof uses \eqn{EqE2b}, \eqn{EqE2c} and
Gauss's contiguity relations for hypergeometric functions
(\cite[\S2.1.2]{Erd}, \cite[\S14.7]{WW}).
There are nine $E_2(i,j)$ terms in \eqn{EqE2d},
and each of them is given by either \eqn{EqE2b} or \eqn{EqE2c},
depending on the relationship between $i$ and $j$.
This means that six separate cases must be considered,
according to whether $k \ge 2n+1$, $k=2n, 2n-1, 2n-2, 2n-3$ or
$k \le 2n-4$.
We give the details just for the first case, the other cases
being very similar.
Assuming then that $k \ge 2n+1$, \eqn{EqE2b} applies to all
nine $E_2(i,j)$ terms in \eqn{EqE2d}.
Writing $\et=k-n$ as before, and replacing the final argument
$\frac{8}{3}$ in the hypergeometric functions by a new
variable $z$, we must show that the expression
\begin{align}\label{EqE2e}
 \frac{(\et+n)!}{(\et-n)! (2n-\et)! 2^{n} 3^{\et-n} } ~ &
  {}_2F_{1}\left[ \begin{array}{c}
                  \et/2-n,\et/2-n+1/2 \\
                       \et-n+1
                  \end{array}
                   ;
                  \begin{array}{c}
                  z
                  \end{array}
                  \right] \nonumber \\
 - ~ \frac{9 n^2 - 9 n + 2}{2} ~ 
\frac{(\et+n-3)!}{(\et-n-1)! (2n-\et)! 2^{n-1} 3^{\et-n-1} } ~ &
   {}_2F_{1}\left[\begin{array}{c}
                  \et/2-n,\et/2-n+1/2 \\
                       \et-n
                  \end{array}
                   ;
                  \begin{array}{c}
                  z
                  \end{array}
                  \right]  \nonumber \\
 + ~ \frac{5}{2} ~
\frac{(\et+n-1)!}{(\et-n+1)! (2n-\et-2)! 2^{n-1} 3^{\et-n+1} } ~ &
   {}_2F_{1}\left[\begin{array}{c}
                  \et/2-n+1,\et/2-n+3/2 \\
                       \et-n+2
                  \end{array}
                   ;
                  \begin{array}{c}
                  z
                  \end{array}
                  \right] \nonumber \\
 - ~ \frac{9 n^2 - 36 n + 35}{2} ~
\frac{(\et+n-4)!}{(\et-n)! (2n-\et-2)! 2^{n-2} 3^{\et-n} } ~ &
   {}_2F_{1}\left[\begin{array}{c}
                  \et/2-n+1,\et/2-n+3/2 \\
                       \et-n+1
                  \end{array}
                   ;
                  \begin{array}{c}
                  z
                  \end{array}
                  \right] \nonumber \\
 - ~ 6(n-1) ~
\frac{(\et+n-3)!}{(\et-n+1)! (2n-\et-3)! 2^{n-2} 3^{\et-n+1} } ~ &
   {}_2F_{1}\left[\begin{array}{c}
                  \et/2-n+3/2,\et/2-n+2 \\
                       \et-n+2
                  \end{array}
                   ;
                  \begin{array}{c}
                  z
                  \end{array}
                  \right] \nonumber \\
 + ~ \frac{3}{2} 
\frac{(\et+n-2)!}{(\et-n+2)! (2n-\et-4)! 2^{n-2} 3^{\et-n+2} } ~ &
   {}_2F_{1}\left[\begin{array}{c}
                  \et/2-n+2,\et/2-n+5/2 \\
                       \et-n+3
                  \end{array}
                   ;
                  \begin{array}{c}
                  z
                  \end{array}
                  \right] \nonumber \\
 - ~ 3(2n-5) 
\frac{(\et+n-4)!}{(\et-n+2)! (2n-\et-5)! 2^{n-3} 3^{\et-n+2} } ~ &
   {}_2F_{1}\left[\begin{array}{c}
                  \et/2-n+5/2,\et/2-n+3 \\
                       \et-n+3
                  \end{array}
                   ;
                  \begin{array}{c}
                  z
                  \end{array}
                  \right] \nonumber \\ 
 - ~ \frac{5}{2}
\frac{(\et+n-3)!}{(\et-n+3)! (2n-\et-6)! 2^{n-3} 3^{\et-n+3} } ~ &
   {}_2F_{1}\left[\begin{array}{c}
                  \et/2-n+3,\et/2-n+7/2 \\
                       \et-n+4
                  \end{array}
                   ;
                  \begin{array}{c}
                  z
                  \end{array}
                  \right] \nonumber \\
 - ~ \frac{5}{2}
\frac{(\et+n-4)!}{(\et-n+4)! (2n-\et-8)! 2^{n-4} 3^{\et-n+4} } ~ &
   {}_2F_{1}\left[\begin{array}{c}
                  \et/2-n+4,\et/2-n+9/2 \\
                       \et-n+5
                  \end{array}
                   ;
                  \begin{array}{c}
                  z
                  \end{array}
                  \right] 
\end{align}
vanishes when $z = \frac{8}{3}$:
Using Gauss's contiguity relations, the nine hypergeometric
functions in \eqn{EqE2e} can all be expressed as linear combinations
of just two of them.
The computer algebra program Maple 11 simplifies\footnote{We don't actually
know how Maple obtains \eqn{EqE2f}, but the result is consistent
with the use of Gauss's relations.}
the above expression to
\begin{align}\label{EqE2f}
 \frac{(\et+n-4)! (3z-8)}
{324 (\et-n+1)! (2n-\et-2)! 2^n 3^{\et-n} z^3 (z-1)^3} 
& \left(
 \phi_1 ~ {}_2F_{1}\left[ \begin{array}{c}
                  \et/2-n+1,\et/2-n+3/2 \\
                       \et-n+2
                  \end{array}
                   ;
                  \begin{array}{c}
                  z
                  \end{array}
                  \right] \right. \nonumber \\
 ~+~ & \left.  \phi_2 ~ {}_2F_{1}\left[ \begin{array}{c}
                  \et/2-n,\et/2-n+1/2 \\
                       \et-n+1
                  \end{array}
                   ;
                  \begin{array}{c}
                  z
                  \end{array}
                  \right] \right) \,,
\end{align}
where $\phi_1$ and $\phi_2$  are polynomials in $z$ of degrees $6$ and $5$
respectively, with coefficients which are polynomials in $n$ and $\et$.
Since the exact values of $\phi_1$ and $\phi_2$ are not important for the 
argument, we relegate them to Tables 
9 and 10
in the Appendix. 
The above expression clearly vanishes 
for $z = \frac{8}{3}$, which proves the desired result.

Second proof. Let 
\beql{EqD2a}
D_2(n,k,c) := \frac{k!} {(2n-k+c)! (k-n-2c)! c! \, 2^{k-n-c} 3^c}
\eeq
denote the first summand in \eqn{EqE2a}.
We look for a recurrence of the form
\beql{EqD2b}
\sum_{r=0}^{4}
\sum_{s=0}^{4}
\sum_{t=0}^{4}
C(r,s,t) D_2(n+r,k+s,c+t) ~=~ 0 \, ,
\eeq
where the coefficients $C(r,s,t)$ depend on $n$ but not 
on $k$ or $c$, with the property that when summed on $c$ it
collapses to the appropriately shifted version of \eqn{EqE2d},
which is:
\begin{align}\label{EqE2dd}
& E_2(n+4,k+4)   - (9 n^2 + 63 n + 110) E_2(n+3,k+1)/2
       + 5 E_2(n+3,k+3)/2 \nonumber \\
& -~ (9 n^2 + 36 n + 35) E_2(n+2,k)/2
       + 6(n+3) E_2(n+2,k+1)
       + 3 E_2(n+2,k+2)/2 \nonumber \\
& -~ 3(2 n+3) E_2(n+1,k)
       - 5 E_2(n+1,k+1)/2
       - 5 E_2(n,k)/2 ~=~ 0 \, .
\end{align}
For this we used the method of Sister Mary Celine Fasenmyer,
exactly as described in \S4.1 of \cite{AeqB}.
A Maple 11 program found that there is a solution to \eqn{EqD2b}
in which the coefficients $C(n,k,c)$ involve five free parameters,
and there is a two-parameter solution which collapses to \eqn{EqE2dd}
when summed on $c$.
The simplest solution (obtained from Maple's solution 
by setting both free parameters to zero) is the following.
All the $C(r,s,t)$ are zero except for the following 19 terms:
\begin{align}
C( 0, 0, 1) &= -8, & C( 2, 1, 1) &= -9, \nonumber \\
C( 0, 0, 2) &= 7, & C( 2, 1, 2) &= 3, \nonumber \\
C( 0, 0, 3) &= -3/2, & C( 2, 2, 1) &= 6, \nonumber \\
C( 1, 0, 1) &= -18, & C( 2, 2, 2) &= -6, \nonumber \\
C( 1, 0, 2) &= 15, & C( 2, 2, 3) &= 3/2, \nonumber \\
C( 1, 0, 3) &= -3, & C( 3, 1, 0) &= -9, \nonumber \\
C( 1, 1, 1) &= -4, & C( 3, 3, 1) &= 5, \nonumber \\
C( 1, 1, 2) &= 3/2, & C( 3, 3, 2) &= -5/2, \nonumber \\
C( 2, 0, 0) &= -9, & C( 4, 4, 1) &= 1, \nonumber \\
C( 2, 0, 1) &= 9. &  \nonumber 
\end{align}
It is easy to verify that this collapses to \eqn{EqE2dd} when
summed on $c$.~~~$\bsq$

\vspace*{+.2in}
Is there a combinatorial proof for \eqn{EqG2d}? We do not know.

We discovered \eqn{EqG2d} by experiment, using Theorem
\ref{th6} to suggest the leading term. (Note that 
if $r(n)$ denotes the right-hand side of \eqn{EqGsf},
then $r(n)/r(n-1) = (9 n^2 - 9 n + 2)/2$.)
We also discovered a second recurrence, which
is independent of \eqn{EqG2d}:
\begin{align}\label{EqG2e}
(n-2) \GG_2(n) & =  n (9 n^2-27 n+17) \GG_2(n-1)/2 \nonumber \\
 & +  (6 n^2-15 n+13/2) \GG_2(n-2) \nonumber \\
 & +  (5 n-5) \GG_2(n-3)/2 \, ,
\end{align}
for $n \ge 3$, with $\GG_2(0)=1$, $\GG_2(1) = 3$,
$\GG_2(2) = 31$.
In view of \eqn{EqG2c}, this is equivalent to a complicated
identity involving hypergeometric functions.
We did not find a proof, but Doron Zeilberger
has kindly informed us that he was able to derive \eqn{EqG2e}
by applying the method of ``creative telescoping''
(\cite[Chap.~6]{AeqB}, \cite{Zeil90b}, \cite{Zeil91})
to \eqn{EqD2a}
and using a modified version of his Maple program ``MultiZeilberger''. 

\section{The case $\si \ge 3$}\label{Sec5}

For $\si \ge 3$ we have not found any connections between
$\GG_{\si}(n)$ and generalized Bessel polynomials or hypergeometric functions,
and we do not have proofs for the recurrences that we have discovered.

However, we do know that recurrences for $\GG_{\si}(n)$ and $E_{\si}(n,k)$ always exist.
This follows from Wilf and Zeilberger's Fundamental Theorem
for Multivariate Sums (\cite[Theorem~4.5.1]{AeqB}, \cite{WZ92a}).

\begin{theorem}\label{th7}
\noindent
(i) For $\si \ge 1$,
there is a number $\delta \ge 0$ such that
$E_{\si}(n,k)$ satisfies a recurrence of the form
\beql{EqEsWZ}
\sum_{i=0}^{\delta} \sum_{j=0}^{\delta} C_{i,j}^{(E)}(n) E_{\si}(n-i,k-j) =0
\mbox{~for~all~} n \,,
\eeq
where the coefficients $C_{i,j}^{(E)}(n)$ are polynomials
in $n$ with coefficients depending on $i$ and $j$.

\noindent
(ii) For $\si \ge 1$,
there is a number $\delta \ge 0$ such that
$\GG_{\si}(n)$ satisfies a recurrence of the form
\beql{EqGsWZ}
\sum_{i=0}^{\delta} \sum_{j=0}^{\delta} C_{i}^{(G)}(n) \GG_{\si}(n-i)=0
\mbox{~for~all~} n \,,
\eeq
where the coefficients $C_{i}^{(G)}(n)$ are polynomials
in $n$ with coefficients depending on $i$.
\end{theorem}

\noindent{\bf Proof.}
(ii) As usual, Eq. \eqn{EqGsWZ} follows by summing \eqn{EqEsWZ} on $k$.
(i) We will use the case $\si = 3$ to illustration of the proof,
the general case being similar. We know from \eqn{EqAAA} that
\beql{EqAAA3}
E_{3}(n,k)
= \sum_{a,b,c,d}
\frac{k!}{a! b! c! d! \, 2^b 6^c 24^d } \,,
\eeq
where the sum is over all $4$-tuples of
nonnegative integers $(a,b,c,d)$ satisfying
\begin{eqnarray}
a + b + c + d & = & n \,, \nonumber \\
a + 2b + 3c + 4d & = & k \,. \nonumber
\end{eqnarray}
In other words,
\beql{EqAAA4}
\frac{E_{3}(n,k)}{2^n}
= \sum_{c,d}
\frac{k!}{(2n-k+c+2d)! (k-n-2c-3d)! c! d! \, 2^{k-c} 3^{c+d} } \,,
\eeq
where now the sum is over all values of $c$ and $d$ for
which the summand is defined.
This summand is a ``holonomic proper-hypergeometric term'',
in the sense of \cite{WZ92a}, and it follows from
the Fundamental Theorem in that paper that 
$E_{3}(n,k)/2^n$ and hence $E_{3}(n,k)$
satisfies a recurrence of the desired form.
Similarly, in the general case, we write the summand in 
$E_{\si}(n,k)$ as a function of $n, k, a_3, \ldots, a_{\si+1}$,
again obtaining a holonomic proper-hypergeometric term.~~~$\bsq$

We conjecture, but do not have a proof,
that a stronger result holds, namely that
recurrences always exist in which the leading terms $C_{0,0}^{(E)}$
and $C_0^{(G)}$ are both $1$, as in \eqn{EqG1d}, \eqn{EqE1d},
\eqn{EqE2d}, \eqn{EqG2d}, \eqn{EqG3e}, \eqn{EqG4e} and 
Table 11.
(The recurrence guaranteed by Theorem \ref{th7} may well look more like
\eqn{EqG2e}, with a nontrivial coefficient on the leading term.)

For $\si = 3, 4$ and $5$, we have found recurrences for
$E_{\si}(n,k)$ and $\GG_{\si}(n)$ with leading coefficient $1$,
although we do not have proofs that they are correct.
The following are our conjectured recurrences for $\GG_3(n)$
and $\GG_4(n)$:

\begin{eqnarray}\label{EqG3e}
 \GG_3(n)  & = & (32 n^3/3 - 16 n^2 + 10 n/3 - 49/6) \GG_3(n-1) \nonumber \\
       & + & (48 n^3 - 236 n^2 + 1157 n/3 - 650/3) \GG_3(n-2) \nonumber \\
       & + & (80 n^3 - 382 n^2 + 641 n - 511) \GG_3(n-3)/3 \nonumber \\
       & + & (64 n^3/3 - 218 n^2 + 2696 n/3 - 7915/6) \GG_3(n-4) \nonumber \\
       & + & (56 n^2 - 490 n + 6853/6) \GG_3(n-5) \nonumber \\
       & + & (56 n - 1703/6) \GG_3(n-6) \nonumber \\
       & + & 58 \GG_3(n-7)/3 \,, 
\end{eqnarray}

\begin{eqnarray}\label{EqG4e}
\GG_4(n) & = &   (625\,{n}^{4}-1250\,{n}^{3}+625\,{n}^{2}-300\,n-543) \GG_4(n-1)/24 \nonumber \\
& + & (27500\,{n}^{4}-184000\,{n}^{3}+447500\,{n}^{2}-473075\,n+180003) \GG_4(n-2) /72 \nonumber \\
& + & (336875\,{n}^{4}-2546500\,{n}^{3}+7679675\,{n}^{2}-12016800\,n+8048577) \GG_4(n-3)/864 \nonumber \\
& + & (4833125\,{n}^{4}-77581625\,{n}^{3}+476892700\,{n}^{2}-1304291160\,n+1325759504) \GG_4(n-4)/2592 \nonumber \\
& + & (1700625\,{n}^{4}+28316750\,{n}^{3}-605973450\,{n}^{2}+3123850885\,n-5033477363) \GG_4(n-5)/7776 \nonumber \\
& + & (2670000\,{n}^{4}-64380500\,{n}^{3}+704577200\,{n}^{2}-3610058445\,n+6818722190)  \GG_4(n-6)/7776 \nonumber \\
& + & (2002500\,{n}^{4}-51976000\,{n}^{3}+517392050\,{n}^{2}-2252744530\,n+3561765885) \GG_4(n-7)/7776 \nonumber \\
& + & (9078000\,{n}^{3}-209915400\,{n}^{2}+1640828980\,n-4301927039) \GG_4(n-8)/7776 \nonumber \\
& + & (5393400\,{n}^{2}-91413680\,n+390747263) \GG_4(n-9)/2592 \nonumber \\
& + & (1593990\,n-14522219) \GG_4(n-11)/972 \nonumber \\
& + & 310343 \GG_4(n-11)/648 \,.
\end{eqnarray}


The recurrence for $\GG_5(n)$ is similar but more complicated,
and we do not state it here. The recurrence
for $E_{3}(n,k)$ is given in the Appendix (see Table 11).
We also omit the recurrences for
for $E_{4}(n,k)$ and $E_{5}(n,k)$, which are even more complicated.

Inspection of these recurrences for
$\si \le 5$ has led us to
some conjectures about their general structure.
First, if $\delta$ denotes the ``depth'' of the recurrence,
as in \eqn{EqEsWZ}, \eqn{EqGsWZ}, then the initial values of $\delta$
for both $\GG_{\si}(n)$ and $E_{\si}(n,k)$
appear to be as shown in Table 2, that is,
it appears that these both recurrences have depth
$\delta = \binom{n+1}{2}+1$ (sequence A000124 of \cite{OEIS}).
\begin{table}[htb]
\label{Tab2}
\caption{Depth $\delta$ of recurrences for $\GG_{\si}(n)$
and $E_{\si}(n,k)$.}
$$
\begin{array}{c|rrrrrr}
\si &    0 & 1 & 2 & 3 & 4 & 5 \\
\hline
\delta & 1 & 2 & 4 & 7 & 11 & 16 
\end{array}
$$
\end{table}
Second, we make the following conjectures\footnote{There are similar
conjectures about the putative recurrence for $\GG_{\si}(n)$.} about
the coefficients in the putative recurrence for $E_{\si}(n,k)$.
We write this recurrence as
\beql{EqEsWZb}
\sum_{i=0}^{\delta} \sum_{j=0}^{\delta} C_{i,j}^{(E)}(n) E_{\si}(n-i,k-j) = 0\,,
\eeq
where $\delta= \binom{n+1}{2}+1, C_{0,0}^{(E)}(n) = 1$.
Then we believe that 
$C_{i,j}(n) = 0$ if $j > \binom{n+1}{2}+1$, or $j<i$, or 
$(i < \si$ and $j > \binom{n+1}{2}+1 - ((\si+1-i)^2 - \si -i-1)/2)$.
Furthermore, the degree of $C_{i,j}(n)$ as a polynomial in $n$ 
is $\le \min \{\si, j-i\}$.

\section{Open questions}\label{Sec6}

We collect here some of the questions that we have mentioned.
(i) The case $\si=1$ corresponds to values of Bessel polynomials; is
there a notion of generalized Bessl polynomial that could be applied
for larger values of $\si$?
(ii) The case $\si=2$ can be described using hypergeometric functions;
is there a notion of generalized hypergeometric function that could be applied
for larger values of $\si$?
(iii) Is there a combinatorial proof of \eqn{EqG2d}?
(iv) Is the conjecture following Theorem \ref{th7} concerning
the existence of recurrences with leading coefficient $1$ true?
(v) Find proofs that the recurrences \eqn{EqG3e} and \eqn{EqG4e}
are correct. 
(vi) Establish the conjectures about the general form of the recurrences for 
$\GG_{\si}(n)$ and $E_{\si}(n,k)$
that are mentioned at the end of \S\ref{Sec5}
(this includes question (iv) as a special case).

\section{Acknowledgment}
We thank Doron Zeilberger for finding a proof
of the recurrence \eqn{EqG2e}.

\noindent{\bf Appendix}

This Appendix collects various tables and multi-line formulas 
that would otherwise have disrupted the flow of the text.

Notation: $E_{\si}(n,k)$ is the number of partitions of $\{1, \ldots, k\}$ into exactly
$n$ blocks of sizes in the range $[1, \ldots, \si+1]$.
Also $\GG_{\si}(n)$ is the number of scenarios when there
are $n+1$ gifts, with a limit of $\si$ steals per gift, and 
the gifts are taken from the pool in the order $1,2,\ldots,n+1$. 

\begin{table}[htbp]
\label{TabE1}
\caption{Values of $E_1(n,k)$.
The array itself appears in several versions in \cite{OEIS}:
see for example A001498, A144299, A144331;
the row sums give A001515, the column sums give A000085 (cf. \cite{MMW}, \cite{MW55}).}
$$
\begin{array}{|c|rrrrrrrrrrrrrrr|} \hline
n \backslash k & 0 & 1 & 2 & 3 & 4 & 5 & 6 & 7 & 8 & 9 & 10 & 11 & 12 & 13 & \ldots  \\
\hline
0 & 1 & 0 & 0 & 0 & 0 & 0 & 0 & 0 & 0 & 0 & 0 & 0 & 0 & 0 & \ldots  \\
1 & 0 & 1 & 1 & 0 & 0 & 0 & 0 & 0 & 0 & 0 & 0 & 0 & 0 & 0 & \ldots  \\
2 & 0 & 0 & 1 & 3 & 3 & 0 & 0 & 0 & 0 & 0 & 0 & 0 & 0 & 0 & \ldots  \\
3 & 0 & 0 & 0 & 1 & 6 & 15 & 15 & 0 & 0 & 0 & 0 & 0 & 0 & 0 & \ldots  \\
4 & 0 & 0 & 0 & 0 & 1 & 10 & 45 & 105 & 105 & 0 & 0 & 0 & 0 & 0 &  \ldots  \\
5 & 0 & 0 & 0 & 0 & 0 & 1 & 15 & 105 & 420 & 945 & 945 & 0 & 0 & 0 & \ldots  \\
6 & 0 & 0 & 0 & 0 & 0 & 0 & 1 & 21 & 210 & 1260 & 4725 & 10395 & 10395 & 0 & \ldots \\
7 & 0 & 0 & 0 & 0 & 0 & 0 & 0 &  1 &  28 &  378 & 3150 & 17325 & 62370 & 135135 & \ldots  \\
8 & 0 & 0 & 0 & 0 & 0 & 0 & 0 & 0 &    1 &  36 &  630 & 6930 & 51975 & 945945 & \ldots  \\
\dots & \dots & \dots & \ldots & \ldots & \ldots & \ldots & \ldots & \ldots & \ldots & \ldots & \ldots & \ldots & \ldots & \ldots & \ldots  \\
\hline
\end{array}
$$
\end{table}

\begin{table}[htbp]
\label{TabE2}
\caption{Values of $E_2(n,k)$.
The array itself appears in several versions in \cite{OEIS}:
see for example A144385, A144399, A144402;
the row sums give A144416, the column sums give A001680 (cf. \cite{MMW}).}
$$
\begin{array}{|c|rrrrrrrrrrrrrrr|} \hline
n \backslash k & 0 & 1 & 2 & 3 & 4 & 5 & 6 & 7 & 8 & 9 & 10 & 11 & 12 & 13 & \ldots  \\
\hline
0 & 1 & 0 & 0 & 0 & 0 & 0 & 0 & 0 & 0 & 0 & 0 & 0 & 0 & 0 & \ldots  \\
1 & 0 & 1 & 1 & 1 & 0 & 0 & 0 & 0 & 0 & 0 & 0 & 0 & 0 & 0 & \ldots  \\
2 & 0 & 0 & 1 & 3 & 7 & 10 & 10 & 0 & 0 & 0 & 0 & 0 & 0 & 0 & \ldots  \\
3 & 0 & 0 & 0 & 1 & 6 & 25 & 75 & 175 & 280 & 280 & 0 & 0 & 0 & 0 & \ldots  \\
4 & 0 & 0 & 0 & 0 & 1 & 10 & 65 & 315 & 1225 & 3780 & 9100 & 15400 & 15400 & 0 & \ldots  \\
5 & 0 & 0 & 0 & 0 & 0 & 1 & 15 & 140 & 980 & 5565 & 26145 & 102025 & 323400 & 800800 & \ldots  \\
6 & 0 & 0 & 0 & 0 & 0 & 0 & 1 & 21 & 266 & 2520 & 19425 & 125895 &  695695 &  3273270 & \ldots  \\
7 & 0 & 0 & 0 & 0 & 0 & 0 & 0 & 1 & 28 & 462 & 5670 & 56595 & 478170 & 3488485 & \ldots  \\
8 & 0 & 0 & 0 & 0 & 0 & 0 & 0 & 0 & 1 & 36 & 750 & 11550 & 144375 & 1531530 & \ldots  \\
\dots & \dots & \dots & \ldots & \ldots & \ldots & \ldots & \ldots & \ldots & \ldots & \ldots & \ldots & \ldots & \ldots & \ldots & \ldots  \\
\hline
\end{array}
$$
\end{table}

\begin{table}[htbp]
\label{TabE3}
\caption{Values of $E_3(n,k)$.
The array itself appears in several versions in \cite{OEIS}:
see A144643, A144644, A144645;
the row sums give A144508, the column sums give A001681 (cf. \cite{MMW}).}
$$
\begin{array}{|c|rrrrrrrrrrrrrrr|} \hline
n \backslash k & 0 & 1 & 2 & 3 & 4 & 5 & 6 & 7 & 8 & 9 & 10 & 11 & 12 & 13 & \ldots  \\
\hline
0 & 1 & 0 & 0 & 0 & 0 & 0 & 0 & 0 & 0 & 0 & 0 & 0 & 0 & 0 & \ldots  \\
1 & 0 & 1 & 1 & 1 & 1 & 0 & 0 & 0 & 0 & 0 & 0 & 0 & 0 & 0 & \ldots \\
2 & 0 & 0 & 1 & 3 & 7 & 15 & 25 & 35 & 35 & 0 & 0 & 0 & 0 & 0 & \ldots \\
3 & 0 & 0 & 0 & 1 & 6 & 25 & 90 & 280 & 770 & 1855 & 3675 & 5775 & 5775 & 0 & \ldots \\
4 & 0 & 0 & 0 & 0 & 1 & 10 & 65 & 350 & 1645 & 6930 & 26425 & 90475 & 275275 & 725725 & \ldots \\
5 & 0 & 0 & 0 & 0 & 0 & 1 & 15 & 140 & 1050 & 6825 & 39795 & 211750 & 1033725 & 4629625 & \ldots \\
6 & 0 & 0 & 0 & 0 & 0 & 0 & 1 & 21 & 266 & 2646 & 22575 & 172095 & 1198120 & 7702695 & \ldots \\
7 & 0 & 0 & 0 & 0 & 0 & 0 & 0 & 1 & 28 & 462 & 5880 & 63525 & 609840 & 5335330 & \ldots \\
8 & 0 & 0 & 0 & 0 & 0 & 0 & 0 & 0 & 1 & 36 & 750 & 11880 & 158235 & 1861860 & \ldots \\
\dots & \dots & \dots & \ldots & \ldots & \ldots & \ldots & \ldots & \ldots & \ldots & \ldots & \ldots & \ldots & \ldots & \ldots & \ldots  \\
\hline
\end{array}
$$
\end{table}

\begin{table}[htbp]
\label{TabE4}
\caption{Values of $E_4(n,k)$.
The array itself appears in several versions in \cite{OEIS}:
see A151338, A151509, A151511;
the row sums give A144509, the column sums give A110038 (cf. \cite{MMW}).}
$$
\begin{array}{|c|rrrrrrrrrrrrrrr|} \hline
n \backslash k & 0 & 1 & 2 & 3 & 4 & 5 & 6 & 7 & 8 & 9 & 10 & 11 & 12 & 13 & \ldots  \\
\hline
0 & 1 & 0 & 0 & 0 & 0 & 0 & 0 & 0 & 0 & 0 & 0 & 0 & 0 & 0 & \ldots  \\
1 & 0 & 1 & 1 & 1 & 1 & 1 & 0 & 0 & 0 & 0 & 0 & 0 & 0 & 0 & \ldots \\
2 & 0 & 0 & 1 & 3 & 7 & 15 & 31 & 56 & 91 & 126 & 126 & 0 & 0 & 0 & \ldots \\
3 & 0 & 0 & 0 & 1 & 6 & 25 & 90 & 301 & 938 & 2737 & 7455 & 18711 & 41811 & 81081 & \ldots \\
4 & 0 & 0 & 0 & 0 & 1 & 10 & 65 & 350 & 1701 & 7686 & 32725 & 132055 & 505351 & 1824823 & \ldots \\
5 & 0 & 0 & 0 & 0 & 0 & 1 & 15 & 140 & 1050 & 6951 & 42315 & 241780 & 1310925 & 6782776 & \ldots \\
6 & 0 & 0 & 0 & 0 & 0 & 0 & 1 & 21 & 266 & 2646 & 22827 & 179025 & 1309000 & 9054045 & \ldots \\
7 & 0 & 0 & 0 & 0 & 0 & 0 & 0 & 1 & 28 & 462 & 5880 & 63987 & 626472 & 5677672 & \ldots \\
8 & 0 & 0 & 0 & 0 & 0 & 0 & 0 & 0 & 1 & 36 & 750 & 11880 & 159027 & 1897896 & \ldots \\
\dots & \dots & \dots & \ldots & \ldots & \ldots & \ldots & \ldots & \ldots & \ldots & \ldots & \ldots & \ldots & \ldots & \ldots & \ldots  \\
\hline
\end{array}
$$
\end{table}

\begin{table}[htbp]
\label{TabE5}
\caption{Values of $E_5(n,k)$.
The array itself appears in several versions in \cite{OEIS}:
see A151359, A151511, A151512;
the row sums give A149187, the column sums give A148092 (cf. \cite{MMW}).}
$$
\begin{array}{|c|rrrrrrrrrrrrrrr|} \hline
n \backslash k & 0 & 1 & 2 & 3 & 4 & 5 & 6 & 7 & 8 & 9 & 10 & 11 & 12 & 13 & \ldots  \\
\hline
0 & 1 & 0 & 0 & 0 & 0 & 0 & 0 & 0 & 0 & 0 & 0 & 0 & 0 & 0 & \ldots  \\
1 & 0 & 1 & 1 & 1 & 1 & 1 & 1 & 0 & 0 & 0 & 0 & 0 & 0 & 0 & \ldots \\
2 & 0 & 0 & 1 & 3 & 7 & 15 & 31 & 63 & 119 & 210 & 336 & 462 & 462 & 0 & \ldots \\
3 & 0 & 0 & 0 & 1 & 6 & 25 & 90 & 301 & 966 & 2989 & 8925 & 25641 & 70455 & 183183 & \ldots \\
4 & 0 & 0 & 0 & 0 & 1 & 10 & 65 & 350 & 1701 & 7770 & 33985 & 143605 & 588511 & 2341339 & \ldots \\
5 & 0 & 0 & 0 & 0 & 0 & 1 & 15 & 140 & 1050 & 6951 & 42525 & 246400 & 1370985 & 7383376 & \ldots \\
6 & 0 & 0 & 0 & 0 & 0 & 0 & 1 & 21 & 266 & 2646 & 22827 & 179487 & 1322860 & 9294285 & \ldots \\
7 & 0 & 0 & 0 & 0 & 0 & 0 & 0 & 1 & 28 & 462 & 5880 & 63987 & 627396 & 5713708 & \ldots \\
8 & 0 & 0 & 0 & 0 & 0 & 0 & 0 & 0 & 1 & 36 & 750 & 11880 & 159027 & 1899612 & \ldots \\
\dots & \dots & \dots & \ldots & \ldots & \ldots & \ldots & \ldots & \ldots & \ldots & \ldots & \ldots & \ldots & \ldots & \ldots & \ldots  \\
\hline
\end{array}
$$
\end{table}

\begin{table}[htbp]
\label{TabG}
\caption{Number $\GG_{\si}(n)$ of scenarios when there are $n+1$ gifts,
a limit of $\si$ steals per gift and 
the gifts are taken from the pool in the order $1,2,\ldots,n$.
(The array itself is entry A144512 in \cite{OEIS};
the rows give A001515, A144416, A144508, A144509, A149187;
the columns A048775, A144511, A144662, A147984.)}
$$
\begin{array}{|c|rrrrrrr|} \hline
n & 0 & 1 & 2 & 3 & 4 & 5 & \ldots  \\
\hline
\GG_0(n) & 1 & 1 & 1 & 1 & 1 & 1 & \ldots  \\
\GG_1(n) & 1 & 2 & 7 & 37 & 266 & 2431 & \ldots  \\
\GG_2(n) & 1 & 3 & 31 & 842 & 45296 & 4061871 & \ldots  \\
\GG_3(n) & 1 & 4 & 121 & 18252 & 7958726 & 7528988476 & \ldots  \\
\GG_4(n) & 1 & 5 & 456 & 405408 & 1495388159 & 15467641899285 & \ldots  \\
\GG_5(n) & 1 & 6 & 1709 & 9268549 & 295887993624 & 34155922905682979 & \ldots \\
\GG_6(n) & 1 & 7 & 6427 & 216864652 & 60790021361170 & 79397199549271412737  & \ldots \\
\GG_7(n) & 1 & 8 & 24301 & 5165454442 & 12845435390707724 & 191739533381111401455478 & \dots \\
\GG_8(n) & 1 & 9 & 92368 & 124762262630 & 2774049143394729653 & 476872353039366288373555323 & \ldots \\
\ldots & \ldots & \ldots & \ldots & \ldots & \ldots & \ldots & \ldots  \\
\hline
\end{array}
$$
\end{table}

%
%
%

\begin{table}
\label{Tabp1}
\caption{The polynomial $\phi_1$ mentioned in the first proof of Theorem \ref{th6}.}
\begin{align}
 \phi_1 & = ( 486+729\,\et n-2349\,n+2916\,{n}^{2}-162\,\et -729\,{n}^{3}-729\,{n}^{2}\et  ) { z}^{6} \nonumber \\
& + \, ( -45\,n+306\,\et -1080\,{n}^{2}+207\,{n}^{3}+180\,{\et }^{3}-216\,{\et }^{2}-324 \,{\et }^{2}n \nonumber \\
& ~~~~~~ -1539\,\et n+1647\,{n}^{2}\et -54 ) {z}^{5} \nonumber \\
& + \, ( 3441\,{n}^{3}-4023\,{n}^{2}\et +24621\,n-348\,{\et }^{3}-16884\,{n}^{2}+1260\,{\et }^{2}n  \nonumber \\ 
& ~~~~~~  -10650\,\et -1800\,{\et }^{2}+ 13203\,\et n-9054 ) {z}^{4} \nonumber \\
& + \, ( 341\,{n}^{3}+948\,{\et }^{2}n+8614\,\et -3359\,n+261 \,{n}^{2}\et +984\,{\et }^{2} \nonumber \\
& ~~~~~~ -6081\,\et n-484\,{\et }^{3}+3270 ) {z}^{3} \nonumber \\
& + \, ( -35572\,n- 36712\,\et n-4092\,{n}^{3}+13512\,{\et }^{2}+2572\,{\et }^{3}+14952 \nonumber \\
& ~~~~~~  +25892\,\et +11164\,{n}^{2}\et - 9244\,{\et }^{2}n+22472\,{n}^{2} ) {z}^{2} \nonumber \\
& + \, ( 11200\,{\et }^{2}n-20160\,{\et }^{2}- 3200\,{\et }^{3}-21120-24320\,{n}^{2}+45760\,\et n \nonumber \\
& ~~~~~~ +4160\,{n}^{3} -38080\,\et +42560\,n-12160\,{ n}^{2}\et  ) z \nonumber \\
& + \,(7680+7680\,{n}^{2}+14080\,\et -15360\,\et n-3840\,{\et }^{2}n+1280\,{\et }^{3} \nonumber \\
& ~~~~~~ - \,14080\,n-1280\,{n}^{3}+7680\,{\et }^{2}+3840\,{n}^{2}\et)  \nonumber 
\end{align}
\end{table}

\begin{table}
\label{Tabp2}
\caption{The polynomial $\phi_2$ mentioned in the first proof of Theorem \ref{th6}.}
\begin{align}
 \phi_2 & = ( 27\,{\et }^{2}+216\,\et n-189\,\et +324+189\,{n}^{2}-675\,n ) {z}^{5} \nonumber \\
& + ( -9
\,{\et }^{2}-9\,{n}^{2}-504\,\et n+495\,\et +9\,n+486 ) {z}^{4} \nonumber \\
& + ( -15\,{\et }^{2}+600
\,\et n-1191\,\et -789\,{n}^{2}-3672+3399\,n ) {z}^{3} \nonumber \\
& + ( -243\,{\et }^{2}+408\,\et n-
555\,\et -303\,{n}^{2}-978+1155\,n ) {z}^{2} \nonumber \\
& + ( 560\,{\et }^{2}-1360\,\et n+3040\,\et 
-3440\,n+720\,{n}^{2}+3840 ) z \nonumber \\
& +(-320\,{\et }^{2}+640\,\et n-1600\,\et -320\,{n}^{2}-1920+
1600\,n) \,. \nonumber
\end{align}
\end{table}

\begin{table}
\label{EqE3e}
\caption{The conjectured recurrence for $E_3(n,k)$. Summing both sides
on $k$ gives \eqn{EqG3e}.}
\begin{align}
E_3(n,k) & = (32 n^3/3-16 n^2+22 n/3-1) E_3(n-1,k-4) \nonumber \\
&        - (4 n+3/2) E_3(n-1,k-2) \nonumber \\
&        - 17 E_3(n-1,k-1)/3 \nonumber \\
&  \nonumber \\
&        + (16 n^3-88 n^2+159 n-189/2) E_3(n-2,k-6) \nonumber \\
&        + (32 n^3-176 n^2+914 n/3-497/3) E_3(n-2,k-5) \nonumber \\
&        + (28 n^2-66 n+46) E_3(n-2,k-4) \nonumber \\
&        + (-12 n+29/2) E_3(n-2,k-3) \nonumber \\
&        - 17 E_3(n-2,k-2) \nonumber \\
&  \nonumber \\
&        + (-16 n^3/3+152 n^2/3-479 n/3+1001/6) E_3(n-3,k-7) \nonumber \\
&        + (32 n^3-262 n^2+2218 n/3-4255/6) E_3(n-3,k-6) \nonumber \\
&        + (84 n^2-382 n+1247/3) E_3(n-3,k-5) \nonumber \\
&        + (16 n-47/3) E_3(n-3,k-4) \nonumber \\
&        - 28 E_3(n-3,k-3) \nonumber \\
&  \nonumber \\
&        + (64 n^3/3-302 n^2+4154 n/3-12427/6) E_3(n-4,k-7) \nonumber \\
&        + (84 n^2-562 n+2858/3) E_3(n-4,k-6) \nonumber \\
&        + (76 n-187) E_3(n-4,k-5) \nonumber \\
&        - 41 E_3(n-4,k-4)/3 \nonumber \\
&  \nonumber \\
&        + (56 n^2-574 n+4352/3) E_3(n-5,k-7) \nonumber \\
&        + (84 n-651/2) E_3(n-5,k-6) \nonumber \\
&        + 17 E_3(n-5,k-5) \nonumber \\
&  \nonumber \\
&        + (56 n-1877/6) E_3(n-6,k-7) \nonumber \\
&        + 29 E_3(n-6,k-6) \nonumber \\
&  \nonumber \\
&        + 58 E_3(n-7,k-7)/3 \,. \nonumber 
\end{align}
\end{table}

\end{document}